\documentclass[titlepage,11pt]{article}

\usepackage[margin=1.2in]{geometry}
\usepackage{latexsym}
\usepackage{amsfonts}
\usepackage[leqno]{amsmath}

\makeatletter
\newcommand{\leqnomode}{\tagsleft@true}
\newcommand{\reqnomode}{\tagsleft@false}
\makeatother

\newcommand{\dichi}{\vec{\chi}}

\usepackage{tikz}
\usepackage{float}
\usepackage{url}
\usetikzlibrary{arrows}
\usetikzlibrary{decorations.markings}
\usepackage[symbol]{footmisc}

\usepackage{xspace}
\usepackage[color=Olive]{todonotes}

\usepackage[colorlinks=true,allcolors=blue]{hyperref}

\newcommand\blackslug{\hbox{\hskip 1pt \vrule width 4pt height 8pt depth 1.5pt
        \hskip 1pt}}
\newcommand\bbox{\hfill \quad \blackslug \bigbreak}

\usepackage{authblk}

\title{Digraphs with all induced directed cycles of the same length are not $\dichi$-bounded}

\author[1]{Alvaro Carbonero}
\author[1]{Patrick Hompe}
\author[2]{Benjamin Moore}
\author[1]{Sophie Spirkl\thanks{Emails: (ar2carbo, phompe, sspirkl)@uwaterloo.ca, brmoore@iuuk.mff.cuni.cz \\
We acknowledge the support of the Natural Sciences and Engineering Research Council of Canada (NSERC), [funding reference number RGPIN-2020-03912]. Cette recherche a été financée par le Conseil de recherches en sciences naturelles et en génie du Canada (CRSNG), [numéro de référence RGPIN-2020-03912].
Benjamin Moore is supported by the ERC-CZ project LL2005 (Algorithms and complexity within and
beyond bounded expansion) of the Ministry of Education of Czech Republic.
}}

\affil[1]{University of Waterloo, Department of Combinatorics and Optimization, Waterloo, Canada}

\affil[2]{Charles University, Institute of Computer Science, Prague, Czech Republic}

\date{\today}

\newtheorem{thm}{Theorem}[section]

\newtheorem{lemma}[thm]{Lemma}

\newcommand{\Proof}{\noindent{\bf Proof.}\ \ }

\begin{document}
\maketitle
\begin{abstract}
For $t \ge 2$, let us call a digraph $D$ \emph{t-chordal} if all induced directed cycles in $D$ have length equal to $t$. In \cite{our_paper}, we asked for which $t$ it is true that $t$-chordal graphs with bounded clique number have bounded dichromatic number. Recently in \cite{chordal}, Aboulker, Bousquet, and de Verclos answered this in the negative for $t=3$, that is, they gave a construction of $3$-chordal digraphs with clique number at most $3$ and arbitrarily large dichromatic number. In this paper,  we extend their result, giving for each $t \ge 3$ a construction of $t$-chordal digraphs with clique number at most $3$ and arbitrarily large dichromatic number, thus answering our question in the negative. On the other hand, we show that a more restricted class, digraphs with no induced directed cycle of length less than $t$, and no induced directed $t$-vertex path, have bounded dichromatic number if their clique number is bounded. We also show the following complexity result: for fixed $t \ge 2$, the problem of determining whether a digraph is $t$-chordal is coNP-complete.

\end{abstract}

\section{Introduction and preliminaries}
Throughout this paper, graphs and digraphs are simple, and in particular, for two vertices $u, v$ in a digraph, not both of the edges $uv$ and $vu$ are present. 

A \emph{$k$-coloring} of a graph $G$ is a function $f:V(G) \to [k]$ such that for each edge $e =xy \in E(G)$, we have $f(x) \neq f(y)$. The \emph{chromatic number} of a graph $G$, denoted $\chi(G)$, is the smallest $k$ such that $G$ admits a $k$-coloring. In this paper we are interested in an analogue of coloring for digraphs called dicoloring first introduced in  \cite{dichi1, dichi2}. 

 We say a digraph is \emph{acyclic} if it has no directed cycle. A \emph{$k$-dicoloring} of a digraph $D$ is a function $f:V(D) \rightarrow [k]$ such that for each $i \in [k]$, the set $f^{-1}(i)$ induces an acyclic subdigraph of $D$. Given a digraph $D$, the \emph{dichromatic number of $D$}, denoted $\dichi(D)$, is the minimum value of $k$ such that there exists a $k$-dicoloring of $D$. 

For a graph $G$, we let $\omega(G)$ be the maximum size of a clique in $G$. For a digraph $D$, we let $\omega(D)$ be the maximum size of a clique in the undirected underlying graph of $D$ (that is, the graph obtained from $D$ with the same vertex set, and two vertices $u$ and $v$ are adjacent if either $uv$ or $vu$ is an edge in $D$). A class $\mathcal{C}$ of graphs is \emph{$\chi$-bounded} (see \cite{gyarfas}) if there exists a function $f$ such that $\chi(G) \leq f(\omega(G))$ for all $G \in \mathcal{C}$.  Similarly, a class $\mathcal{C}$ of digraphs is \emph{$\dichi$-bounded} if there exists a function $f$ such that $\dichi(D) \le f(\omega(D))$ for all $D \in \mathcal{C}$. 

In \cite{our_paper}, motivated by problems for the dichromatic number, we proved that the class of graphs where each triangle-free induced subgraph has bounded chromatic number is not $\chi$-bounded. Since we were interested in the dichromatic number we posed various open problems. In particular, we asked if the class of digraphs where every induced directed cycle has length $t$, for a fixed integer $t \geq 3$, is $\dichi$-bounded. We will refer to a digraph in which all induced directed cycles have length $t$ as \emph{$t$-chordal}. 

In \cite{chordal}, Aboulker, Bousquet, and de Verclos answered our question in the negative for $t=3$. Our main result is to extend this negative answer to all $t \geq 3$:

\begin{thm}\label{main1}
For $t \ge 3$, and every $N \in \mathbb{N}$, there exists a $t$-chordal digraph $D$ with $\omega(D) \le 3$ and $\dichi(D) \ge N$, and if $t > 3$, then $\omega(D) \leq 2$.
\end{thm}

We give a positive result, which in some sense demonstrates that Theorem \ref{main1} is tight, by proving $\dichi$-boundedness for a subclass of $t$-chordal digraphs:

\begin{thm}\label{our_result}
For every integer $t$, there is a function $f_t: \mathbb{N} \rightarrow \mathbb{N}$ such that for every digraph $D$ which is $t$-chordal and has no induced directed path with exactly $t$ vertices, we have $\dichi(D) \leq f_t(\omega(D))$.
\end{thm}

Since directed cycles of length more than $t$ contain a directed path on $t$ vertices, Theorem \ref{our_result} is equivalent to saying that for every $t$, digraphs with no induced $t$-vertex path and no induced directed cycle on strictly less than $t$ vertices are $\dichi$-bounded. 

Despite not being $\dichi$-bounded, one could hope that $t$-chordal digraphs still have a ``nice" structural description. We show that any clean structural description of $t$-chordal digraphs is unlikely, as deciding if a digraph is $t$-chordal is coNP-complete:

\begin{thm}\label{main2}
Fix an integer $t \geq 2$. Deciding if a given digraph is $t$-chordal is coNP-complete. 
\end{thm}

The paper is organized as follows. In Section \ref{mainsec}, we prove Theorem \ref{main1} and Theorem \ref{our_result}. In Section \ref{hardnesssec}, we prove Theorem \ref{main2}.

\section{The construction}

Our construction uses a key idea similar to and inspired by \cite{chordal}, making sure that in every $k$-dicoloring, certain independent sets miss a color, and arranging this for one independent set at a time. This is accomplished through the following lemma, from which Theorem \ref{main1} will be derived. 
\label{mainsec}

\begin{lemma}\label{key_thm}
For $t \ge 3$, suppose that $D$ is a $t$-chordal digraph with $\omega(D) \le 3$ and $\dichi(D) = k$. If $\mathcal{C} = \{I_1,\ldots,I_p\}$ is a collection of independent sets in $D$, then there exists a $t$-chordal digraph $D'$ with the following properties:
\begin{itemize}
    \item $\omega(D') \le 3$, and if $\omega(D) \leq 2$ and $t > 3$, then $\omega(D') \leq 2$; and
    \item for every $k$-dicoloring of $D'$, there exists a copy of $D$ as an induced subgraph of $D'$ such that for each $1 \le i \le p$, the copy of $I_i$ contained in that copy of $D$ is dicolored with at most $k-1$ colors.
\end{itemize} 
\end{lemma}
\Proof We will prove the lemma by induction on the chromatic number of an auxiliary graph that we define now. For a digraph $D$ and a collection $\mathcal{C} = \{I_1, \dots, I_p\}$ of independent sets, let $G_{D,\mathcal{C}}$ be the graph with $V(G_{D,\mathcal{C}}) = \{1,2,\ldots,p\}$, and where $ij \in E(G)$ if and only if $I_{i} \cap I_{j} \ne \emptyset$; in other words, $G_{D,\mathcal{C}}$ is the intersection graph of $\mathcal{C}$. 

The base case is when $\chi(G_{D,\mathcal{C}}) = 0$, (and thus, $\mathcal{C} = \emptyset$) where the statement is trivially true. Now, for $s \ge 0$, suppose that the statement is true for all digraphs $D$ paired with a collection of independent sets $\mathcal{C}$ such that the corresponding graph $G_{D,\mathcal{C}}$ has $\chi(G_{D,\mathcal{C}}) \le s$, and suppose that $\chi(G_{D,\mathcal{C}}) = s+1$. Take an $(s+1)$-coloring of $G_{D,\mathcal{C}}$, say $f$. Let $X_{0} = \{j_{1},\ldots,j_{q}\}$ be a color class of $f$, and let $S_0 = \{I_{j_1},\ldots,I_{j_q}\}$. It follows that $\chi(G_{D,\mathcal{C} \setminus S_0}) = s$, and that the sets $I_{j_1}, \dots, I_{j_q}$ are pairwise disjoint. 

We define a sequence $D_0, \dots, D_q$ with the following properties:
\begin{itemize}
    \item $D_{0} = D$;
    \item for all $i \in \{0, \dots, q\}$, $D_i$ is a digraph with clique number at most 3 (at most 2 if $\omega(D) \leq 2$ and $t > 3$);  
    \item for all $i \in \{0, \dots, q\}$, $D_i$ contains pairwise disjoint copies $G^1_i, \dots, G^{t^i}_i$ of $D$;
    \item for all $i \in \{0, \dots, q\}$, in every $k$-dicoloring of $D_i$, there is a $t^* \in \{1, \dots, t^i\}$ such that for every $r \in \{1, \dots, i\}$, the independent set corresponding to $I_{j_r}$ in the copy $G^{t^*}_i$ of $D$ is dicolored with at most $k-1$ colors; and
    \item for all $i \in \{0, \dots, q\}$ and for every $r \in \{i+1, \dots, q\}$, the union of the independent sets corresponding to $I_{j_r}$ in $G_i^1, \dots, G_i^{t^i}$ is an independent set in $D_i$. 
\end{itemize}

Clearly, $D_0$ satisfies the properties above for $i = 0$. Suppose that we have defined $D_i$ for some $i \in \{0,\ldots,q-1\}$.  Let $I$ be the union of the sets corresponding to $I_{j_{i+1}}$ in $G_i^1, \dots, G_i^{t^i}$. From the properties of $D_i$, it follows that $I$ is an independent set. We create a new digraph $D_{i+1}$ as follows. Let  $D_{i}^{1},\ldots,D_{i}^{t}$ be $t$ copies of $D_{i}$, and let $V(D_{i+1}) = \bigcup_{j=1}^{t} V(D_{i}^{j})$. In addition to the edges corresponding to each copy $D_i^j$ of $D_i$, we add the following edges. For $j \in \{1,\ldots,t\}$, let $I^j$ denote the copy of $I$ in $D_i^j$. Then, for each vertex $v \in I^{j}$ and $u \in I^{j+1}$, we add the edge $vu$ (where indices are taken modulo $t$, so $I^{t+1}$ means $I^1$). 

Observe that the clique number of $D_{i+1}$ is at most $3$, and at most $2$ if $t > 3$ and $\omega(D) \leq 2$. This follows since a clique in $D_{i+1}$ either has all its vertices in $D^{j}_{i}$ for some $j$, or it has at most one vertex contained in each copy of $D^{j}_{i}$. Similarly, $D_{i+1}$ is $t$-chordal, as every induced directed cycle either uses precisely one vertex from each $D^{j}_{i}$, or it is contained completely in a copy $D^{j}_{i}$ for some $j$. 

Since $D_{i+1}$ contains $t$ copies of $D_i$, it follows that $D_{i+1}$ contains $t$ times as many copies of $D$ as $D_i$ does, for $t \cdot t^i = t^{i+1}$ copies overall. Let us label them as $G_{i+1}^1, \dots, G_{i+1}^{t^{i+1}}$ arbitrarily. It follows that the last bullet holds for $D_{i+1}$, since it holds for $D_i$, and since the only edges between copies of $D_i$ in $D_{i+1}$ are between vertices that are in a set corresponding to $I_{j_{i+1}}$ in a copy of $D$, and $I_{j_{i+2}}, \dots, I_{j_q}$ are disjoint from $I_{j_{i+1}}$ from the choice of $S_0$. 

Now, let $c$ be a $k$-dicoloring of $D_{i+1}$ (if one exists). Then, not all $k$ colors occur in each of $I^1, \dots, I^t$, for otherwise there would be a monochromatic $t$-vertex directed cycle. It follows that there is a $j \in \{1, \dots, t\}$ such that the copy of $I$ in $D_i^j$ is dicolored with at most $k-1$ colors. Consequently, for each of the copies of $G_i^1, \dots, G_i^{t^i}$ of $D$ in $D_i^j$, the copy of $I_{j_i}$ in this copy of $D$ is dicolored with at most $k-1$ colors. Together with the fact that the third bullet holds for $D_i$, by applying it to $D_i^j$, it follows that it holds for $D_{i+1}$. 

This completes the definition of $D_0, \dots, D_q$. Note that $D_q$ has the property that in every $k$-dicoloring of $D_q$, there is a copy of $D$ in $D_q$ such that each of $I_{j_1}, \dots, I_{j_q}$ is dicolored with at most $k-1$ colors. 

Now we construct a collection of independent sets $\mathcal{C}'$. For every independent set $I \in \mathcal{C} \setminus S_0$, and for every $t^* \in \{1, \dots, t^q\}$, we add the copy of $I$ in $G_q^{t^*}$ to $\mathcal{C}'$. Note that we can $s$-color $G_{D_q, \mathcal{C}'}$ by fixing an $s$-coloring of $G_{D, \mathcal{C}\setminus S_0}$, and assigning to each independent set in a copy of $D$ the color of the independent set in $D$ it corresponds to. This, by construction, is an $s$-coloring, since two sets which are assigned the same color are either copies of disjoint sets, or disjoint copies of the same set.

It follows that we can apply the inductive hypothesis to $D_q$ and $\mathcal{C}'$; and so there is a digraph $D^*$ with $\omega(D^*) \leq 3$ (and if $\omega(D) \leq 2$ and $t > 3$ then $\omega(D^*) \leq 2$), and for every $k$-dicoloring of $D^*$, there is a copy $D_q'$ of $D_q$ in $D^*$ such that for every $I \in \mathcal{C}'$, the copy of $I$ in $D_q'$ is dicolored with at most $k-1$ colors. It follows that in $D_q'$, there are copies $G_q^1, \dots, G_q^{t^q}$ of $D$, and for every $I \in \mathcal{C} \setminus S_0$, every copy of $I$ in every one of $G_q^1, \dots, G_q^{t^q}$ receives at most $k-1$ colors. But also, from the properties of $D_q$, there is a $t^* \in \{1, \dots, t^q\}$ such that for $G_q^{t^*}$ in $D_q'$, for each $I \in S_0$, the copy of $I$ in $G_q^{t^*}$ receives at most $k-1$ colors. This shows that in this copy  $G_q^{t^*}$ in $D_q'$ of $D$, the second bullet of the lemma holds. This concludes the proof. \bbox{}

Now, we present our main result, which is proved in the same way as in \cite{chordal} (in which the main theorem  is derived from their Lemma 3 exactly like we derive ours from Lemma \ref{key_thm} above):
\begin{thm}
For all $t \ge 3$, there exists a sequence of $t$-chordal digraphs $\{D_n\}$ such that for all $n \ge 1$, we have $\omega(D_n) \le 3$ (and $\omega(D_n) \leq 2$ if $t > 3$) and $\dichi(D_n) \ge n$.
\end{thm}
\Proof Let $D_1$ be the digraph with one vertex and no edges, and define the sequence of digraphs inductively as follows. For $k \ge 1$, take $\dichi(D_k)$ disjoint copies of $D_k$, forming a digraph $D_k'$, and construct a collection $\mathcal{C}$ of independent sets of $D_k'$ as follows: for each set of $\dichi(D_k)$ vertices, one in each copy of $D_k$ in $D_k'$, place this set in $\mathcal{C}$. Then by Lemma \ref{key_thm}, we have that there exists a $t$-chordal digraph $D_{k+1}$ with $\omega(D_{k+1}) \le 3$ (and $\omega(D_n) \leq 2$ if $t > 3$) such that for any $\dichi(D_k)$-dicoloring of $D_{k+1}$ there exists a copy of $D_k'$ such that each independent set in $\mathcal{C}$ uses at most $\dichi(D_{k})-1$ colors. Now, for a $\dichi(D_k)$-dicoloring of $D_k'$, it follows that in each of the $\dichi(D_k)$ copies of $D_k$ in $D_k'$ there is a vertex of each color, so there exists an independent set in $\mathcal{C}$ with exactly one vertex of each color, which is a contradiction. It follows that there does not exist a $\dichi(D_k)$-dicoloring of $D_{k+1}$ and thus $\dichi(D_{k+1}) \ge \dichi(D_k)+1$, and therefore $\dichi(D_k) \ge k$ for all $k \ge 1$, as desired. This completes the inductive step and completes the proof.\bbox{}

We now show the following positive result which shows that the above construction is in some sense tight:

\begin{thm} \label{thm:end}
For every $l$, there is a function $f_l: \mathbb{N} \rightarrow \mathbb{N}$ such that for every digraph $D$ with no induced directed cycle of length less than $l$ and no induced directed path with exactly $l$ vertices, we have $\dichi(D) \leq f_l(\omega(D))$. 
\end{thm}
\Proof
Let us denote by $\mathcal{C}_l$ the class of digraphs with no induced directed cycle of length less than $l$ and no induced directed path with exactly $l$ vertices.

We will show that $f_l(\omega) = (l+1)^\omega$ gives the desired result. We proceed by induction on $k = \omega(D)$. For $k=1$, the statement is true, since in that case $\dichi(D) = 1 \le l+1 = (l+1)^\omega$.

Now, suppose that the statement is true for $k = n$, namely that for all digraphs $D' \in \mathcal{C}_l$ with $\omega(D') = n$, we have $\dichi(D') \le (l+1)^n$. Now, consider a digraph $D \in \mathcal{C}_l$ with $\omega(D) = n + 1$. We will show that $\dichi(D) \le (l+1)^{n+1}$. Suppose for the sake of a contradiction that $\dichi(D) > (l+1)^{n+1}$. First, note that we may assume that $D$ is strongly connected, since the dichromatic number of a digraph is equal to the maximum of the dichromatic numbers of its strongly connected components. Let $v_1 \in V(D)$ be an arbitrary vertex, and let $C_1,\ldots,C_s$ be the strongly connected components of $D \setminus ( \{v_1\} \cup N(v_1))$, where $N(v_1)$ denotes the set of vertices which have an in-edge or an out-edge to $v_1$. We claim that there exists a component $C_i$ with $\dichi(C_i) > (l+1)^{n+1} - (l+1)^n$. Indeed, since $\omega(N(v_1)) \le n$, we have by induction that $\dichi(N(v_1)) \le (l+1)^n$. So if $\dichi(C_i) \le (l+1)^{n+1} - (l+1)^n$ for all $1 \le i \le t$, we can dicolor each component $C_i$ with the same set of $(l+1)^{n+1}-(l+1)^n$ colors, dicolor $N(v_1)$ with a disjoint set of $(l+1)^n$ colors, and reuse one of the colors from the $C_i$ to dicolor $v_1$. Clearly, this is a dicoloring of $D$ with at most $(l+1)^{n+1}$ colors, which is a contradiction. Thus we may assume without loss of generality that $\dichi(C_1) > (l+1)^{n+1} - (l+1)^n$.

Now, since $D$ is strongly connected, there exists a directed path from $v_1$ to a vertex in $C_1$. Let $P'$ be the shortest such path, and let the second-to-last vertex in $P'$ be $v_2$. Now, we let $P$ be the portion of the path $P'$ from $v_1$ to $v_2$. First, we note that $P$ is forward-induced, meaning that there are no edges in the direction of the path other than those in $P$. Furthermore, since $P'$ was the shortest directed path from $v_1$ to $C_1$, it follows that $v_2$ is the only vertex in $P$ which has an out-edge to a vertex in $C_1$.

Let $D_1 = D$. We define a sequence of digraphs as follows. Let $D_2 = D[\{v_2\} \cup C_1]$, and note that since $C_1$ is strongly connected and $v_2$ has an out-edge to $C_1$, it follows that there exists a directed path from $v_2$ to every vertex in $D_2$. Then we iterate this argument, letting $C_1'$ be a strongly connected component of $D_2 \setminus (\{v_2\} \cup N(v_2))$ with $\dichi(C_1') > (l+1)^{n+1} - 2 (l+1)^n$, taking a shortest directed path $Q$ from $v_2$ to $C_1'$, and adding to our path $P$ the subpath of $Q$ from $v_2$ to the second-to-last vertex of $Q$, which we call $v_3$. Since $(l+1)^{n+1}- (l+1) (l+1)^n \geq 0$, we can iterate this argument $l+1$ times and obtain a forward-induced path $P$ with a subset of the vertices of $P$ specified as $v_1,\ldots,v_{l+2}$. In particular, the $l$ vertices $v_3,\ldots,v_{l+2}$ are not in $\{v_1\} \cup N(v_1)$,  since all vertices of $P$ after $v_2$ are in $C_1$, and thus not in $\{v_1\} \cup N(v_1)$.

Now, let $w$ be the last vertex in $P$ which is in $N(v_1)$. Then it follows from the properties of $P$ that there exists a forward-induced path $P'$ with $l$ vertices starting at $w$ such that $w$ is the only vertex in $P'$ which is in $\{v_1\} \cup N(v_1)$. Since $P'$ is a forward-induced path on $l$ vertices and there is no induced directed path on $l$ vertices and no induced directed cycle of length less than $l$, it follows that $P'$ is an induced directed cycle of length $l$. Let the vertices of $P'$ be $(w_1,\ldots,w_l)$ such that $w_1=w$. Now, if $w_1$ is an out-neighbor of $v_1$, then we have that $(v_1 w_1 w_2, \ldots, w_{l-1})$ is an induced directed path on $l$ vertices, which is a contradiction. Suppose instead that $w_1$ is an in-neighbor of $v_1$. Then we have that $(w_3 w_4, \ldots, w_l w_1 v_1)$ is an induced directed path on $l$ vertices, which is a contradiction as well. Hence, $\dichi(D) \le (l+1)^n$. This finishes the proof by induction and shows that $\dichi(D) \le (l+1)^{\omega(D)}$ for all digraphs $D \in \mathcal{C}_l$, as desired. \bbox{}

\section{Deciding if a graph is $t$-chordal is coNP-complete}
\label{hardnesssec}
In this section, we prove Theorem \ref{main2}. The construction is very similar to those of \cite{bienstock, lubiw, kk}; in particular, \cite{lubiw} showed that for every $r \geq 0$ and $d \geq 2$, deciding if a given digraph has an induced directed cycle of residue $r$ modulo $d$ is NP-complete. Note that Theorem \ref{thm:nphard} can be modified in a straightforward way to show that the problem of determining if a non-simple digraph (where both $uv$ and $vu$ are allowed to be present at the same time) is $2$-chordal is coNP-complete as well, whereas we can decide in polynomial time if a simple digraph is $2$-chordal by checking if it is acyclic. 
\begin{thm} \label{thm:nphard}
For $t \ge 3$, the problem of determining whether a digraph $D$ is $t$-chordal is coNP-complete.
\end{thm}
\Proof Observe the problem is in coNP. To see this, given any digraph $D$, a certificate that it is not $t$-chordal is an induced directed cycle $C$ which has length not equal to $t$. Further, such a certificate can be checked in polynomial time.

For the rest of the proof, let $t \geq 3$ be a fixed integer. We will give a reduction from $3$-SAT, a well-known NP-complete problem. Let $\phi$ be an instance of $3$-SAT, with variables $x_{1},\ldots,x_{n}$ and clauses $C_{1},\ldots,C_{m}$. We will construct a digraph $D$ such that $D$ is $t$-chordal if and only if $\phi$ is a NO-instance of $3$-SAT, which suffices to prove the result. 

We build $D$ in stages. First, for each variable $x_{i}$ $i \in \{1,\ldots,n\}$, we create a variable gadget, $D_{x_{i}}$, which is a digraph with two vertices $v^{i}_{1},v^{i}_{2}$ and two directed paths $P^{i}_{1}$, $P^{i}_{2}$, both directed from $v^{i}_{1}$ to $v^{i}_{2}$, which are vertex-disjoint aside from the vertices $v^{i}_{j}$, $j \in \{1,2\}$ and both contain exactly $t$ edges. For each $j \in \{1,2\}$, let $z^{i}_{j}$ be the vertex in $P^{i}_{j}$ adjacent from $v^{i}_{1}$ and $q^{i}_{j}$ the vertex adjacent to $v^{i}_{2}$ in $P^{i}_{j}$. 

Second, for each clause $C_{i}$, $i \in \{1,\ldots,m\}$, we create a clause gadget, $D_{C_{i}}$ which has two vertices $u^{i}_{1},u^{i}_{2}$, and for each literal $y$ contained in $C_{i}$, we add a directed path $u^{i}_{1},w^{y},u^{i}_{2}$. 

Now we create a digraph $D$ whose vertex set is $$\{V(D_{x_{i}}) : i \in \{1,\ldots,n\} \} \cup \{V(D_{C_{i}}) : i \in \{1,\ldots,m\}\}.$$ For the edges,  for $i \in \{1,\ldots,n-1\}$, we add the edge $v^{i}_{2}v^{i+1}_{1}$, as well as the edge $v^{n}_{2}u^{1}_{1}$. For $j \in \{1,\ldots,m-1\}$ we add the edges $u^{i}_{2}u^{i+1}_{1}$, as well as the edge $u^{m}_{2}v^{1}_{1}$. For a clause $C_{i}$, if the non-negated variable $x_{i}$ appears in the clause, then we add edges $w^{x_{i}}z^{i}_{2}$ and $q^{i}_{2}w^{x_{i}}$.  For a clause $C_{i}$ if the negated variable $x_{i}$ appears in the clause, then we add edges $w^{x_{i}}z^{i}_{1}$ and $q^{i}_{2}w^{x_{i}}$. 

Now we claim that $D$ is $t$-chordal if and only if the $3$-SAT instance has no solution. Suppose there is an induced directed cycle $C$ of length other than $t$. Suppose first that $C$ uses an edge of the form $w^{x_{i}}z^{i}_{j}$ or $q^{i}_{j}w^{x_{i}}$ for some $j \in \{1, 2\}$. If $C$ uses $w^{x_{i}}z^{i}_{j}$, then, since $z^{i}_{j}$ has a unique out-neighbor, it follows that $C$ contains the path $P^i_j \setminus v_2^i$. In particular, $C$ contains $q^i_j$; and since $C$ contains $w^{x_i}$, and $C$ is induced, it follows that $C$ contains $q^i_jw^{x_i}$, and so $C$ has length $t$. If $q^{i}_{j}w^{x_{i}}$, then, since $q_j^i$ has a unique in-neighbor, it follows similarly that $C$ has length $t$. Now, by construction, it follows $C$ has non empty-intersection with every variable gadget and every clause gadget.

For each variable $x_{i}$, if $C$ takes the path $P^{i}_{1}$, we set $x_{i}$ to true, and otherwise we set $x_{i}$ to false. We claim this gives a satisfying assignment to the $3$-SAT instance. If not, then some clause evaluates to false, say $C_{1}$ without loss of generality. We may assume up to symmetry that $C \cap D_{C_{1}}$ is the path associated to the variable $x_{i}$. If this variable is not negated in $C$ then this implies that $C$ uses the path $P^{i}_{2}$ to set $x_{i}$ to false, however in this case there are edges from vertices in $P^{i}_{2}$ to $w^{x_{i}}$ which are not apart of $C$, contradicting that $C$ is induced. An analogous argument holds when $x_{i}$ is negated in $C_{1}$. 

Similarly for the converse, if an assignment to $x_{1},\ldots,x_{m}$ satisfies the formula, we construct the induced directed cycle by adding either $P^{i}_{1}$ if $x_{i}$ is true, or $P^{i}_{2}$ if $x_{i}$ is false, taking any path from each clause gadget which corresponds to a true literal in that clause, and adding the edges between clause gadgets and variable gadgets in the natural way. It is easy to check that this directed cycle is an induced directed cycle of length not equal to $t$, completing the proof.\bbox{}  

\section*{Acknowledgments}
We acknowledge the support of the Natural Sciences and Engineering Research Council of Canada (NSERC), [funding reference number RGPIN-2020-03912]. Cette recherche a été financée par le Conseil de recherches en sciences naturelles et en génie du Canada (CRSNG), [numéro de référence RGPIN-2020-03912].

\end{document}